\documentclass{amsart}
\usepackage{amssymb}
\usepackage{a4wide}
\usepackage{enumerate}
\usepackage{graphics}
\title{Covers of Point-Hyperplane Graphs}
\author{Arjeh~M.~Cohen}
\author{E.J.~Postma}

\newtheorem{theorem}{Theorem}[section]
\newtheorem{lemma}[theorem]{Lemma}
\newtheorem{corollary}[theorem]{Corollary}

\providecommand{\myinput}[1]{\includegraphics{#1.eps}}

\providecommand{\card}[1]{\left\lvert#1\right\rvert}
\providecommand{\class}[1]{\widetilde{#1}}
\providecommand{\bigw}{{\textstyle\bigwedge}}
\providecommand{\tensor}{\otimes}
\providecommand{\rightcosets}[2]{#1\setminus\!\!\setminus\,#2}

\DeclareMathOperator{\characteristic}{char}
\DeclareMathOperator{\Aut}{Aut}

\begin{document}
\begin{abstract}
  We construct a cover of the non-incident point-hyperplane graph of
  projective dimension 3 for fields of characteristic 2. If the
  cardinality of the field is larger than 2, we obtain an elementary
  construction of the non-split extension of $\mathrm{SL}_4
  (\mathbb{F})$ by $\mathbb{F}^6$.
\end{abstract}
\maketitle

\section{Introduction}

The non-incident point-hyperplane graph $H_n (\mathbb{F})$ has as
vertex set the non-incident pairs of a point and a hyperplane in the
projective geometry of projective dimension $n$ over a field
$\mathbb{F}$. Two distinct vertices are adjacent if the points and
hyperplanes are mutually incident. These graphs have been studied
extensively, cf.~Gramlich~\cite{cit:gramlich}. One of their properties
is that $H_{n + 1} (\mathbb{F})$ is locally $H_n (\mathbb{F})$ for all
$\mathbb{F}$ and $n$, and that every connected and locally $H_n
(\mathbb{F})$ graph is isomorphic to $H_{n + 1} (\mathbb{F})$ whenever
$n > 2$.

This property does not necessarily hold if $n \leq 2$. Indeed, if $n
\in \{0, 1\}$ then it is easily seen not to hold. In
Gramlich~\cite{cit:gramlich}, a covering graph of $H_3(\mathbb{F}_2)$,
constructed by means of a computer algebra computation, shows that it
does not hold for $n = 2, \mathbb{F} = \mathbb{F}_2$. In this paper we
give a computer-free construction of a covering graph of $H_3
(\mathbb{F})$ for $\characteristic \mathbb{F} = 2$, thus providing
counterexamples to the local recognizability of $H_3 (\mathbb{F})$ for
this wider class. This is the content of our main theorem,
Theorem~\ref{thm:mainthm}, which is based on a construction developed
in Sections~\ref{s:labelcover} and~\ref{s:phpgraphs}. These sections are
based on the second author's Master's thesis~\cite{cit:msc}.

In Section~\ref{s:autgroup}, we find automorphisms of this covering
graph, generating an extension of $\textrm{SL}_4 (\mathbb{F})$ by
$\mathbb{F}^6$, which is non-split if $\card{\mathbb{F}} > 2$.
Consequently, as a byproduct we find an elementary construction of the
non-split extension discussed by~Bell~\cite{cit:bell} and
Griess~\cite{cit:griess}.

\subsection{Notation and conventions}
We let groups act on the right. Hence we will often consider the set
of right cosets of a subgroup $H$ of some group $G$. In order to avoid
confusion with the set minus operation, we will denote this set as 
$\rightcosets{H}{G}$.

All graphs are simple and undirected. The adjacency of vertices $u$
and $v$ is denoted by $u \perp v$. For a graph $\Gamma$, we let $V
(\Gamma)$ be the set of its vertices and $D (\Gamma)$ be the set of
its \emph{darts} or \emph{oriented} edges; that is, the set of ordered
pairs of vertices $(u, v)$ for which $u \perp v$.

\section{The voltage assignment} \label{s:labelcover}

In this section, we discuss a general method of constructing covers of
a given graph by means of voltage assignments. For a general
introduction to voltage assignments, see Malni\v
c~et~al.~\cite{cit:malnic}. 

For all vertices $v$ of a graph, we call the induced graph on the
neighbourhood of $v$ the \emph{local graph} at $v$. Let $\Gamma$ and
$\Delta$ be two graphs. A map $\alpha\colon V(\Gamma) \to
V(\Delta)$ preserving adjacency such that the local graph at every
vertex of $\Gamma$ is mapped isomorphically to the local graph at its
image, is called a \emph{local isomorphism}. If $\Gamma$ and $\Delta$
are connected and a local isomorphism from $\Gamma$ to $\Delta$
exists, we call $\Gamma$ a \emph{cover} of $\Delta$.

Let $N$ be a group. A map $\ell\colon D (\Delta) \to N$ such that
$\ell (u, v) = \ell (v, u)^{-1}$ is called a \emph{voltage assignment}
of $\Delta$. We will often write $\ell_{u, v}$ for $\ell (u, v)$, or
$\ell_{i, j}$ for $\ell_{v_i, v_j}$. The \emph{lift} of $\Delta$ with
respect to $\ell$ is the graph with vertex set $V (\Delta) \times N$,
where $(u, m) \perp (v, n)$ if and only if $u \perp v$ and $\ell (u,
v) = m n^{-1}$.

Given a path $P = (v_0 \perp v_1 \perp \dotsb \perp v_i)$, we call
$\ell_{0, 1} \ell_{1, 2} \dotsm \ell_{n - 1, n}$ the \emph{voltage} of
$P$, denoted by $\ell (P)$.  Using induction it is immediate that for
any $m \in N$, we find an induced path in the lift from $(v_0,
m)$ to $(v_n, \ell (P)^{-1} m)$.

Let $\Delta$ be a connected graph with voltage assignment $\ell\colon
D(\Delta) \to N$. Let $\Gamma$ be the lift of $\Delta$ with respect to
$\ell$. Then it is easy to see that $\Gamma$ is connected if and only
if for every $n \in N$ and every $v_0 \in V (\Delta)$, there is an $i
\in \mathbb{N}$ and a cycle $(v_0, v_1, \ldots, v_i = v_0)$ such that
$\ell_{0,1} \ell_{1,2} \cdots \ell_{i - 1, i} = n$. It is equally easy
to see that for all $n \in N$ and $v \in V(\Delta)$, the local graph
at $(v, n)$ in $\Gamma$ is isomorphic to the local graph at $v$ in
$\Delta$, if and only if for every triangle $u, v, w$, we have
$\ell_{u, v} \ell_{v, w} \ell_{w, u} = 1$.  These two observations
lead to the following straightforward lemma.

\begin{lemma} \label{lm:sublabelingconditions}
  Let $\Delta$ be a connected graph with voltage assignment $\ell'\colon
  D(\Delta) \to N'$. Let $T$ be the normal closure of the subgroup of
  $N'$ generated by the voltages of all triangles.
  
  Define $N = \rightcosets{T}{N'}$ and $\ell_{u, v} = T \ell_{u, v}'$.
  Let $M$ be the subgroup of $N$ generated by the voltages (with
  respect to $\ell$) of all cycles. Let $\Gamma$ be the lift of
  $\Delta$ with respect to $\ell$.
  
  Then the map $\alpha \colon V(\Gamma) \to V(\Delta),\ (v, n) \mapsto
  v$ is a local isomorphism and every connected component of $\Gamma$
  is an $\card{M}$-fold cover of $\Delta$.
\end{lemma}
Let $G$ be a group of automorphisms of $\Delta$ with an action on $N$.
We will say that $\ell$ is \emph{$G$-equivariant} if and only if for
all $g \in G$ and $v \perp w \in \Delta$, we have that $\ell_{v^g,
  w^g} = (\ell_{v, w})^g$.

Group-equivariant voltage assignments enable the group to lift to a
group of automorphisms of the lift. This is the content of the next
lemma, the proof of which is again straightforward.
\begin{lemma} \label{lm:labelledaction}
  Let $G$ be a subgroup of $\Aut \Delta$ such that $\ell$ is
  $G$-equivariant. Let $\Gamma$ be the lift of $\Delta$ with
  respect to $\ell$. Then $G \ltimes N$ acts faithfully on $\Gamma$ by
  the action $(v, n)^{(g, k)} = (v^g, n^g k)$.
\end{lemma}
Now suppose we have the setup of Lemma~\ref{lm:sublabelingconditions}.
Suppose that $M$ is Abelian and that $\ell$ is $G$-equivariant. Choose
a vertex $v \in V(\Delta)$.  For all $g \in G$, choose a path $P_g$
from $v^g$ to $v$, and let $\lambda (g)$ be the voltage of $P_g$. Then
the following lemma holds.
\begin{lemma} \label{lm:GM}
  The stabilizer in $G \ltimes N$ of the connected component
  $\Gamma_0$ of $\Gamma$ containing $(v, 0)$ is $H = \{(g, \lambda (g)
  + m) \mid g \in G, m \in M\}$, which is an extension of $G$ by $M$.
\end{lemma}
\begin{proof}
  Since $(v, 0)^{(g, \lambda (g) + m)} = (v^g, \lambda (g) + m)$ and
  since the path induced on $P_g$ starting at $(v^g, \lambda (g) + m)$
  ends at $(v, m)$, we have that $H$ stabilizes $\Gamma_0$. If any
  element $(g, n)$ stabilizes $\Gamma_0$, it maps $(v, 0)$ to an
  element $(v^g, n)$ such that there is a path from $(v^g,
  \lambda(g))$ to $(v^g, n)$. Then the projection of that path down to
  $\Delta$ is a cycle; hence $\lambda - n \in M$. So $H$ is the full
  stabilizer of $\Gamma_0$. Therefore it is a group.

  The kernel of the projection onto the first coordinate is $\{1\}
  \times M$, so that is a normal subgroup. The quotient by that
  subgroup is $G$.
\end{proof}
\subsection{The reduct} \label{ss:reduct}
Let a graph $\widetilde \Delta$ be given. We define an equivalence
relation $\sim$ on the vertices by
\[ v \sim w \qquad \Longleftrightarrow \qquad v^\perp = w^\perp, \]
where $v^\perp$ is the set of neighbours of $v$. The adjacency
relation carries over in a natural way from the vertices to the
equivalence classes of $\sim$ (which are necessarily cocliques). Hence
we can mod out $\sim$ and obtain a new graph $\widetilde \Delta / \sim
\ = \Delta$. We call $\Delta$ the \emph{reduct} of $\widetilde
\Delta$, and write $\class v$ for the vertex in $\Delta$ representing
the equivalence class containing $v$.

We will need this well-known construction in Section~\ref{s:phpgraphs}
in order to link the projective graph $H_n (\mathbb{F})$ to its affine
version, $\widetilde H_n (\mathbb{F})$ (defined in the same section).

\begin{lemma} \label{lm:tildesim}
  Let $\ell\colon D (\widetilde \Delta) \to N$ be a voltage assignment
  on a graph $\widetilde \Delta$ without isolated vertices. Let
  $\widetilde \Gamma$ be the lift of $\widetilde \Delta$ with respect
  to $\ell$. Then the following assertions are equivalent:
  \begin{enumerate}[(i)]
  \item \label{it:lisl} for all $u, v, w \in V(\widetilde \Delta)$
    such that $u \sim v$ and $v \perp w$, we have $\ell_{w, u} =
    \ell_{w, v}$,
  \item \label{it:equiv} for all $(u, m), (v, n) \in V(\widetilde
    \Gamma)$, we have that $(u, m) \sim (v, n)$ if and only if $m = n$
    and $u \sim v$.
  \end{enumerate}
\end{lemma}
\begin{proof}
  (\textit{\ref{it:equiv}}$\,$) $\Rightarrow$ (\textit{\ref{it:lisl}}$\,$)\ :
  Let $u$, $v$, $w$ be such that $u \sim v \perp w$. Then
  $(u, 1) \sim (v, 1)$.  Now $(u, 1) \perp (w, \ell_{w, u})$ and $(v,
  1) \perp (w, \ell_{w, v})$; this implies that $(u, 1) \perp (w,
  \ell_{w, v})$.  But $(u, 1)$ has only one neighbour of shape $(w,
  -)$. Hence $\ell_{w, u} = \ell_{w, v}$.
  
  (\textit{\ref{it:lisl}}$\,$) $\Rightarrow$ (\textit{\ref{it:equiv}}$\,$)\ :
  Suppose that $u \sim v \perp w$ for $u, v, w \in V(\widetilde
  \Delta)$. Then for any $n \in N$, both $(u, n)$ and $(v, n)$ are
  adjacent to $(w, \ell_{w,u} n)$. Similarly for other neighbours of
  $u$ and $v$. So $(u, n) \sim (v, n)$.
  
  Now let us assume $(u, m) \sim (v, n)$. Then $u \sim v$.  Since
  $\widetilde \Delta$ contains no isolated vertices, there exists a
  $(w, j)$ such that $(v, n) \perp (w, j)$; then $\ell_{v, w} = n
  j^{-1}$, and also $\ell_{u,w} = m j^{-1}$.  Since the two are the
  same, we find that indeed $m = n$.
\end{proof}

If the conditions of Lemma~\ref{lm:tildesim} hold, we will say that
$\ell$ is \emph{reductive}. If $\ell$ is reductive, then by
(\textit{\ref{it:equiv}}$\,$) we can write $(\class{v}, n)$ for
$\class{(v, n)}$; and then by (\textit{\ref{it:lisl}}$\,$),
$\ell_{\class u, \class v}$ is well defined as $\ell_{u, v}$.

In the following three corollaries, $\widetilde \Delta$ will be a
connected graph. $\ell\colon D(\Delta) \to N$ will be a reductive
voltage assignment. $\widetilde \Gamma$ will be the lift of
$\widetilde \Delta$ with respect to $\ell$.  Furthermore, $\Delta =
\widetilde \Delta / \sim$ and $\Gamma = \widetilde \Gamma / \sim$.

\begin{corollary}
  Define $ \ell'\colon D (\Delta) \to N$ by $\ell_{\class u, \class
    v}' = \ell_{u, v}$. Then $\Gamma$ is the lift of $\Delta$ with
  respect to $\ell'$.
\end{corollary}
\begin{corollary} \label{cor:tildelabelledaction}
  Let $G$ be a group of automorphisms of $\widetilde \Delta$. Suppose
  $\ell$ is $G$-equivariant. Then $G \ltimes N$ has an action on
  $\Gamma$ defined by $(\class{v}, n)^{(g, k)} = (\class{v}^g, n^g
  k)$.
\end{corollary}
\begin{corollary} \label{cor:subtildelabelingconditions}
  Let $\ell'\colon D (\Delta) \to N'$ be a reductive voltage
  assignment on $\widetilde \Delta$ that is $G$-equivariant for some
  group $G$ of automorphisms of $\widetilde \Delta$. Let $T$ be the
  normal closure of the subgroup of $N'$ generated by the voltages of
  triangles.
  
  Define $N = \rightcosets{T}{N'}$ and $\ell_{u, v} = T\ell_{u, v}'$.
  Let $M$ be the subgroup of $N$ generated by the voltages of cycles.
  Then the map $\alpha \colon \Gamma \to \Delta, (\class{v}, n)
  \mapsto \class{v}$ is a local isomorphism and every connected
  component of $\Gamma$ is an $\card{M}$-fold cover of $\Delta$.
  
  Choose a vertex $v \in \Delta$ and for all $g \in G$, choose a path
  $P_g$ from $v^g$ to $v$, and set $\lambda (g) = \ell (P_g)$. If $M$
  is Abelian, then the stabilizer in $G \ltimes N$ of the connected
  component of $\Gamma$ containing $(v, 0)$ is $\{(g, \lambda (g) + m)
  \mid g \in G, m \in M\}$, which is an extension of $G$ by $M$.
\end{corollary}

\section{Point-Hyperplane Graphs} \label{s:phpgraphs}
Consider the projective geometry $\mathbb{P}_n (\mathbb{F})$ of
(projective) dimension $n$ over the field $\mathbb{F}$. We denote
incidence by $\subset$ or $\supset$ and the projective dimension by
$\dim$. Now define the graph $H_n (\mathbb{F})$ to have
vertex set
\[ 
\{ (x, X) \mid x, X \in \mathbb{P}_n (\mathbb{F}), \dim x = 0, \dim
X = n - 1, x \not \subset X \} 
\] 
and adjacency defined by
\[ 
(x, X) \perp (y, Y) \qquad \Longleftrightarrow \qquad x \subset Y
\quad \textrm{and} \quad y \subset X.
\]
For this graph, with $n = 3$ and $\mathbb{F}$ of characteristic $2$,
we will build a cover as follows. We first construct an affine version
$\widetilde H_3 (\mathbb{F})$, the reduct of which is $H_3
(\mathbb{F})$. Then we recall some multilinear algebra in order to
define a voltage assignment for $\widetilde H_3 (\mathbb{F})$. This
provides an $\card{\mathbb{F}^6}$-fold cover of $\widetilde H_3
(\mathbb{F})$, the reduct of which is an $\card{\mathbb{F}^6}$-fold
cover of $H_3 (\mathbb{F})$.

$V$ will be the vector space $\mathbb{F}^{n + 1}$ with basis $e_1,
e_2, \dotsc$ and dual basis $f_1, f_2, \dotsc$, so $\mathbb{P}_n
(\mathbb{F}) = \mathbb{P} (V)$. We define $\widetilde H_n
(\mathbb{F})$ to be the graph with vertex set
\[ \{ v \tensor f \mid v \in V,\ f \in V^*,\ f (v) \not= 0\}, \]
and adjacency defined by
\[
v \tensor f \ \perp \ w \tensor g \qquad \Longleftrightarrow \qquad f (w)
= v (g) = 0.
\]
\begin{lemma} \label{lm:DeltaisH}
  Let $n \geq 2$. Then $\widetilde H_n (\mathbb{F}) / \sim$ is isomorphic
  to $H_n(\mathbb{F})$.
\end{lemma}
\begin{proof}
  Clearly, for all $\alpha \in \mathbb{F}$, we have $v \tensor f \ 
  \sim \ \alpha (v \tensor f)$.  Also, if $v \tensor f \ \sim \ w
  \tensor g$, then $w$ is in the intersection of all $n$-dimensional
  subspaces of $V$ containing $v$, and $g$ is in the intersection of
  all $n$-dimensional subspaces of $V^*$ containing $f$. Hence $w
  \tensor g = \alpha(v \tensor f)$, for some $\alpha \in
  \mathbb{F}^*$. So the vertex sets of $\widetilde H_n (\mathbb{F})$
  and $H_n (q)$ are equal. Clearly the adjacency relation is also the
  same.
\end{proof}
We will proceed to recall some basic multilinear algebra which will be
needed for constructing the voltage assignment of $\widetilde H_n
(\mathbb{F})$.

We fix $n = 3$ so that $\dim V = 4$ and let
\[
\bigw^k V = V^{\tensor k} / \langle v_1 \tensor v_2 \tensor \dotsb
\tensor v_k \mid v_i = v_j \text{ for some }i \not= j \rangle
\]
be the $k$th Grassmannian of $V$. The image of $v_1 \tensor \dotsb \tensor
v_k$ in $\bigw^k V$ is denoted $v_1 \wedge \dotsb \wedge v_k$. Let $G$
be a group with a linear action on $V$; this induces a natural action
on $\bigw^k V$. Now $G \leq \mathrm{SL} (V)$ if and only if $G$
stabilizes every element of $\bigw^4 V$. We will mostly be using the
case where $k = 2$. We need the following elementary lemmas.

\begin{lemma}
  There is a canonical isomorphism between ${\bigl(\bigw^2 V\bigr)}^*$
  and $\bigw^2 \left(V^*\right)$, that preserves the induced action of
  $\mathrm{GL} (V)$.
\end{lemma}
\begin{proof}[Sketch of proof]
  Let $B \colon \bigw^2 \left(V^*\right) \times
  \bigw^2 V \to \mathbb{F}$ be defined for $\hat f = f_1 \wedge f_2 \in
  \bigw^2 \left(V^*\right)$ and $\hat v = v_1 \wedge v_2 \in \bigw^2
  V$ by
  \[
  B (\hat f, \hat v) = f_1 (v_1) f_2 (v_2) - f_1 (v_2) f_2 (v_1),
  \]
  and extended by linearity. We define $b_{\hat f} (\hat v) = B (\hat
  f, \hat v)$. Then $\beta\colon \hat f \mapsto b_{\hat f}$ is an
  isomorphism.
\end{proof}

Because of the preceding lemma, we will drop the parentheses in the
future and write $\bigw^2 V^*$.

\begin{lemma}
  Let $V$ be a vector space of dimension $4$ over a field
  $\mathbb{F}$. Fix an isomorphism $\chi$ between $\bigw^4 V$ and
  $\mathbb{F}$. Let $G$ be a subgroup of $\mathrm{SL} (V)$. Then there is a
  canonical isomorphism between $\bigw^2 V^*$ and $\bigw^2 V$, which
  respects the natural induced group actions of $G$ on $\bigw^2 V^*$
  and $\bigw^2 V$.
\end{lemma}
\begin{proof}[Sketch of proof]
  We can set up a bilinear mapping $B$ from $\bigw^2 V \times \bigw^2
  V$ to $\mathbb{F}$ as follows:
  \[
  B (v_1 \wedge v_2, w_1 \wedge w_2)
  = (v_1 \wedge v_2 \wedge w_1 \wedge w_2) ^ \chi,
  \]
  and extended by linearity. Now let $\psi$ map $\hat w \in
  \bigw^2 V$ to the linear functional that maps $\hat v \in
  \bigw^2 V$ to $B (\hat v, \hat w)$. Then $\psi$ is a
  vector space isomorphism. 
\end{proof}

We will often let $\phi$ denote the inverse of $\psi$, and if we
consider it appropriate, leave the map out entirely.

For an arbitrary vector space $V$, we let
\[ 
S_2 (V) = (V \tensor V) / \langle v \tensor w - w \tensor v \mid v, w
\in V \rangle
\]
be the second order symmetric tensor of $V$. Then the natural action
of $G$ on $V \tensor V$ induces a natural action on $S_2 (V)$. We
denote the image of $v \tensor w$ in $S_2 (V)$ by $vw$. We will often
write $w^2$ for $ww$.

Now let $\characteristic \mathbb{F} = 2$,
and let $W = \bigw^2 V = \bigw^2 (V^*)$ of dimension $6$. Then $S_2
(W)$ has dimension $21$.

\begin{lemma}
  The subspace $W^{(2)}$ of $S_2 (W)$, defined as
  \[ W^{(2)} = \langle \hat w^2 \mid \hat w \in W \rangle, \]
  has dimension $6$ and is invariant under the induced action of
  $\mathrm{GL} (V)$.
\end{lemma}
\begin{proof}[Sketch of proof]
  If $(\hat e_i)_{i = 1}^6$ is a basis for $W$, then $(\hat e_i^2)_{i
    = 1}^6$ is a basis of $W^{(2)}$. For any $g \in \mathrm{GL} (V)$, $\hat w
  \in W$ we have $(\hat w^2)^g = (\hat w^g)^2$.
\end{proof}
We now fix an isomorphism between $\bigw^4 V$ and $\mathbb{F}$.
\begin{lemma}
  Let $w, x, y, z \in V$ such that $w \wedge x \wedge y \wedge z = 1$.
  Then
  \[ 
  U = (w \wedge x)(y \wedge z) + (w \wedge y)(z \wedge x) + (w \wedge
  z)(x \wedge y)
  \]
  does not depend on the choice of $w, x, y, z$. Furthermore, it is
  fixed under the induced action of $\mathrm{SL} (V)$.
\end{lemma}
\begin{proof}  
  The map
  \[
  \Delta \colon (w, x, y, z) \mapsto (w \wedge x)(y \wedge z) + (w \wedge y)(z
  \wedge x) + (w \wedge z)(x \wedge y)
  \]
  is $4$-linear and alternating. There is only one such map: the
  determinant. Hence for tuples of vectors such that $\det (w, x, y,
  z) = w \wedge x \wedge y \wedge z = 1$, we find that $\Delta$ must
  be constant.
  
  Since the image of $\Delta (w, x, y, z)$ under an element of
  $\mathrm{SL} (V)$ is $\Delta (w', x', y', z')$ for some tuple
  satisfying $w' \wedge x' \wedge y' \wedge z' = 1$, the element $U$
  is fixed under $\mathrm{SL} (V)$.
\end{proof}

\subsection{The voltage assignment $\ell$}
We now restrict our attention to $H_3 (\mathbb{F})$, for some
$\mathbb{F}$ with $\characteristic \mathbb{F} = 2$. Hence $\dim V = 4$
and we retain $W$, $W^{(2)}$ and $U$ as in the previous section.  We
let $\ell\colon D (\widetilde H_3 (\mathbb{F})) \to S_2 (W)$ assign
the voltage
\begin{equation} \label{eq:label}
h_1 (v_1)^{-1} h_2 (v_2)^{-1} (v_1 \wedge v_2)(h_1 \wedge
h_2)^\phi
\end{equation}
to the dart from $v_1 \tensor h_1$ to $v_2 \tensor h_2$, and let
$\ell^U \colon D(\widetilde H_3 (\mathbb{F})) \to S_2 (W) / \langle U
\rangle_{\mathbb{F}_2}$ be the composition of $\ell$ with the natural
projection of $S_2 (W)$ to $S_2 (W) / \langle U \rangle_{\mathbb{F}_2}$.

\begin{theorem} \label{thm:mainthm}
  Let $\widetilde \Gamma$ be the lift of $\widetilde H_3 (\mathbb{F})$
  with respect to $\ell^U$, and let $\Gamma = \widetilde \Gamma /
  \sim$. Then every connected component of $\Gamma$ is an
  $\card{\mathbb{F}^6}$-fold cover of $H_3 (\mathbb{F})$.
\end{theorem}
For proving Theorem~\ref{thm:mainthm}, we need some auxiliary lemmas.
\begin{lemma}
  The voltage assignment $\ell$, defined in~\eqref{eq:label}, is reductive.
\end{lemma}
\begin{proof}
  Let $\widetilde u_1 = v_1 \tensor h_1, \widetilde u_2 = v_2 \tensor
  h_2, \widetilde u_3 = v_3 \tensor h_3$ be vertices of $\widetilde
  H_3 (\mathbb{F})$ satisfying $\widetilde u_1 \sim \widetilde u_2
  \perp \widetilde u_3$. Since $\widetilde u_1 \sim \widetilde u_2$,
  we can write $v_2 \tensor h_2$ as $\alpha v_1 \tensor h_1$ for some
  $\alpha \in \mathbb{F}$.  Then
  \[ 
  \ell_{\widetilde u_2, \widetilde u_3} = \alpha^{-1} h_1 (v_1)^{-1}
  h_3 (v_3)^{-1} \alpha (v_1 \wedge v_3)(h_1 \wedge h_3)^\phi =
  \ell_{\widetilde u_1, \widetilde u_3}.
  \]
\end{proof}
\begin{lemma} \label{lm:triangleinU}
  Each triangle in $\widetilde H_3 (\mathbb{F})$ has voltage $U$.
\end{lemma}
\begin{proof}
  Let $(v_1 \tensor h_1)$, $(v_2 \tensor h_2)$, $(v_3 \tensor h_3)$ be
  a triangle in $\widetilde H_3 (\mathbb{F})$. Then $\{v_i\}$ and
  $\{h_i\}$ are both linearly independent sets. Hence the intersection
  of the null spaces of $\{h_i\}$ has dimension $1$; choose $t$ nonzero
  in it. Then $\{v_1, v_2, v_3, t\}$ is linearly independent. We define
  $u = (v_1 \wedge v_2 \wedge v_3 \wedge t)^{-1} t$.
  Then $\{v_1, v_2, v_3, u\}$ is a basis for $V$.
  
  Since $h_1 (v_3)$, $h_2 (v_3)$, $h_1 (u)$, $h_2 (u)$ all vanish,
  $(h_1 \wedge h_2)^\phi = \alpha (v_3 \wedge u)$ for some nonzero
  $\alpha \in \mathbb{F}$. From this we can derive that 
  $1 = v_1 \wedge v_2 \wedge v_3 \wedge u = \alpha^{-1} \bigl[h_1
  (v_1) h_2 (v_2) + h_1 (v_2) h_2 (v_1) \bigr] = \alpha^{-1} h_1 (v_1)
  h_2 (v_2)$.
  Hence
  $(h_1 \wedge h_2)^\phi = h_1 (v_1) h_2 (v_2) (v_3 \wedge u)$;
  similarly we obtain
  $(h_1 \wedge h_3)^\phi = h_1 (v_1) h_3 (v_3) (v_2 \wedge h_2)$ and
  $(h_2 \wedge h_3)^\phi = h_2 (v_2) h_3 (v_3) (v_1 \wedge h_1)$.
  So if $\{i, j, k\} = \{1, 2, 3\}$, then the voltage of the dart from
  $v_i \tensor h_i$ to $v_j \tensor h_j$ is
  $ h_i (v_i)^{-1} h_j (v_j)^{-1} (v_i \wedge v_j)(h_i \wedge
  h_j)^\phi = (v_i \wedge v_j)(v_k \wedge u)$.
  The sum of the voltages is then
  \[
  (v_1 \wedge v_2)(v_3 \wedge u) + (v_1 \wedge v_3)(v_2 \wedge u) +
  (v_2 \wedge v_3)(v_1 \wedge u) = U.
  \]
\end{proof}
For the following four lemmas, note that we regard $S_2 (W)$ as a
group only here, so the subgroups are the subspaces over
$\mathbb{F}_2$ -- not necessarily over $\mathbb{F}$.
\begin{lemma} \label{lm:quadrangleinWplusU}
  The voltage of a cycle in $\widetilde H_3 (\mathbb{F})$ of length
  four is in the $\mathbb{F}_2$-space $W^{(2)} \oplus \langle U \rangle$.
\end{lemma}
\begin{proof}
  Note that $U \not\in W^{(2)}$, so the space is indeed a direct sum.
  
  Consider a cycle consisting of $v_0 \tensor h_0, \dotsc, v_3 \tensor
  h_3$. We define $\alpha_i = h_i (v_i)^{-1}$. We first assume that
  $v_0 = v_2$. Then the voltage of the quadrangle is equal to $ (v_0
  \wedge v_1) (\alpha_1 h_1 \wedge (\alpha_0 h_0 + \alpha_2 h_2))^\phi
  +$ $(v_0 \wedge v_3) (\alpha_3 h_3 \wedge (\alpha_0 h_0 + \alpha_2
  h_2))^\phi$.  Note that $(\alpha_0 h_0 + \alpha_2 h_2) (v_i) = 0$
  for all $i$: for $i = 1$ or $3$ because of the adjacencies in the
  graph, and for $i = 0$ or $2$ because then $ (\alpha_0 h_0 +
  \alpha_2 h_2) (v_i) = 0$.  Hence the null space of $\alpha_0 h_0 +
  \alpha_2 h_2$ contains all $v_i$.
  
  We will show that we may assume that $(h_i)$ is a linearly
  independent set. Let us assume that a nontrivial linear combination
  of $(h_i)$ is zero. We will distinguish cases according to the sets
  of nonzero coefficients.

  \newlength{\mylength} \setlength{\mylength}{\textwidth}
  \addtolength{\mylength}{-\parindent}
  \vspace{.3\baselineskip}
%  \vspace{1pt}
  \begin{minipage}{\mylength}
    If only the coefficients for two adjacent vertices are nonzero,
    say for $h_0$ and $h_1$, then their null spaces coincide.  Hence
    $v_0$ is in the null space of $h_0$. Contradiction.
%  \end{minipage}\vspace{9pt}
    
%  \begin{minipage}{\mylength}
    If only the coefficients for $h_0$ and $h_2$ are nonzero, then we
    have $v_0 \tensor h_0 \sim v_2 \tensor h_2$; then the voltage is
    the same as if we replace $h_2$ by $h_0$, and if we do that we get
    the sum of two voltages of $2$-cycles, which is $0$.
%  \end{minipage}\vspace{9pt}
    
%  \begin{minipage}{\mylength}
    If only the coefficients for $h_1$ and $h_3$ are nonzero, then say
    $h_3 = \gamma h_1$. The voltage is then $ (v_0 \wedge (\alpha_1
    v_1 + \gamma \alpha_3 v_3))(h_1 \wedge (\alpha_0 h_0 + \alpha_2
    h_2))^\phi$.  But since both $h_1$ and $\alpha_0 h_0 + \alpha_2
    h_2$ vanish on both $v_0$ and $\alpha_1 v_1 + \gamma \alpha_3
    v_3$, the elements $(v_0 \wedge (\alpha_1 v_1 + \gamma \alpha_3
    v_3))$ and $(h_1 \wedge (\alpha_0 h_0 + \alpha_2 h_2))^\phi$ only
    differ by a scalar factor. Hence the voltage is in $W^{(2)}$.
%  \end{minipage}\vspace{9pt}
    
%  \begin{minipage}{\mylength}
    If exactly three of the coefficients are nonzero, say those for
    $h_0$, $h_1$ and $h_2$, then the common null space of any pair of
    those is equal to the common null space of the three. In
    particular, the common null space of $h_0$ and $h_2$ is contained
    in the null space of $h_1$. So $v_1$ is in the null space of
    $h_1$.  Contradiction.
%  \end{minipage}\vspace{9pt}
    
%  \begin{minipage}{\mylength}
    If all four coefficients are nonzero, the common null space of
    $h_1$ and $\alpha_0 h_0 + \alpha_2 h_2$ coincides with the common
    null space of $h_3$ and $\alpha_0 h_0 + \alpha_2 h_2$; hence for
    some nonzero $\lambda, \mu \in \mathbb{F}$ we have $ \alpha_3 h_3
    = \lambda \alpha_1 h_1 + \mu (\alpha_0 h_0 + \alpha_2 h_2)$.  Then
    the voltage of the quadrangle is $ (v_0 \wedge (v_1 + \lambda
    v_3))(\alpha_1 h_1 \wedge (\alpha_0 h_0 + \alpha_2 h_2))^\phi$.
    Again, both $\alpha_1 h_1$ and $\alpha_0 h_0 + \alpha_2 h_2$ are
    $0$ on both $v_0$ and $v_1 + \lambda v_3$. Therefore the voltage
    is in $W^{(2)}$.
  \end{minipage}
  \vspace{0cm}
  
  So $(h_i)$ is a linearly independent set. Then $(\alpha_i h_i)$ is a
  basis for $V^*$. Let $(x_i)$ be its dual basis (so, a basis of $V$).
  We write $x_{i,j}$ for $x_i \wedge x_j$.  We will now determine the
  isomorphism between $\bigw^2 V^*$ and $\bigw^2 V$ with respect to
  these bases. Let $D = (x_0 \wedge \dotsb \wedge x_3)^{-1}$. Let
  $w_0$, $w_1$ be arbitrary vectors in $V$ and let us write $w_j =
  \sum_{i = 0}^3 \eta_{i,j} x_i$.  Now let us see what element of
  $\bigw^2 V^*$, expressed in $(\alpha_i h_i)_{i = 0}^3$, corresponds
  to e.g.~$x_{2,3}$. For clarity, we will use $[$ and $]$ to
  denote application of a functional in this computation.
  \begin{multline} \label{eq:dual}
    (x_2 \wedge x_3) [w_0 \wedge w_1] = (\eta_{0,0} \eta_{1,1} +
    \eta_{0,1} \eta_{1,0}) (x_0 \wedge x_1 \wedge x_2 \wedge x_3) = \\
    D^{-1} (\alpha_0 h_0 [w_0] \alpha_1 h_1 [w_1] + \alpha_0
    h_0 [w_1] \alpha_1 h_1 [w_0]) = D^{-1} (\alpha_0 h_0
    \wedge \alpha_1 h_1) [w_0 \wedge w_1].
  \end{multline}
  We see that $\alpha_i h_i \wedge \alpha_j h_j$ corresponds to
  $D x_{k,\ell}$, where $\{i, j, k, \ell\} = \{1, 2, 3,
  4\}$.

  Let us write
  $v_j = \sum_{i = 0}^3 \epsilon_{i,j} x_i$.
  Then 
  $\epsilon_{i,j} = \alpha_i h_i (v_j) = h_i (v_i)^{-1} h_i (v_j)$;
  in particular, $\epsilon_{i,i} = 1$ and $\epsilon_{i,i \pm 1} =
  0$. Using the fact that $v_0 = v_2$, we can find the values of
  $\epsilon_{i,j}$ to be as in Table~\ref{tab:epsilon}, where $\beta$
  and $\gamma$ are arbitrary scalars.

  \begin{table}
    \centering  
    \caption{The values of $\epsilon_{i,j}$. The values of $\beta$ and
      $\gamma$ are arbitrary.}
    \label{tab:epsilon}
    \[
    \begin{array}{l|*{4}{r}}
      &\multicolumn{4}{c}{i} \\
      j  & 0 & 1       & 2 & 3 \\
      \hline
      0, 2 & 1 & 0       & 1 & 0 \\
      1    & 0 & 1       & 0 & \gamma \\
      3    & 0 & \beta & 0 & 1
    \end{array}
    \]
  \end{table}  
  
  If we now define $x_4 = x_0 + x_2$, then the voltage of the
  quadrangle is equal to
  \[
  (x_4 \wedge (x_1 + \gamma x_3))(x_{2,3} + x_{0,3}) + (x_4 \wedge
  (\beta x_1 + x_3))(x_{1,2} + x_{0,1}), \] which evaluates to $\gamma
  (x_{3,4})^2 + \beta (x_{1,4})^2 \in W^{(2)}$. This proves the lemma
  for quadrangles with $v_0 = v_2$. Duality gives us that the same
  holds if $h_0 = h_2$. Let us call quadrangles for which two opposite
  vertices have the vector or the functional in common,
  \emph{special}.
  
  Now let us assume that we have a non-special quadrangle $v_0 \tensor
  h_0, \dotsc, v_3 \tensor h_3$. If $v_0 \tensor h_2$ is a vertex
  (i.e.~$h_2 (v_0) \not= 0$), then we can split the quadrangle into
  two special quadrangles. See Figure~\ref{fig:splitquadrangle}.
  
  \begin{figure}
    \centering
    \myinput{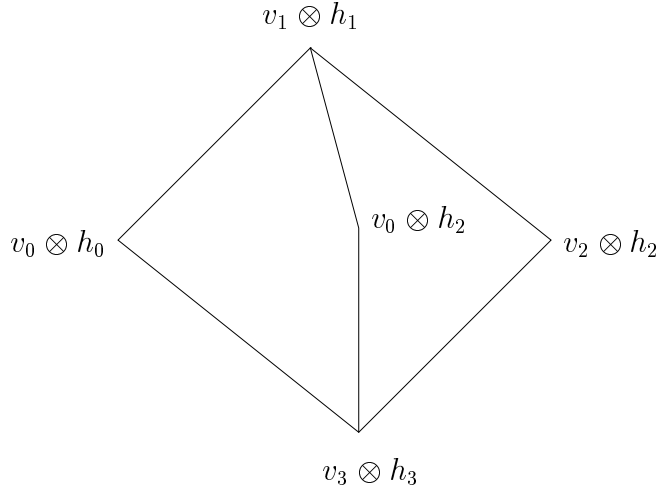}
    \caption{We split a general quadrangle into two special quadrangles.}
    \label{fig:splitquadrangle}
  \end{figure}
  
  This leaves us with the case where for all $i$, we have $h_{i + 2}
  (v_i) = 0$. Then $\alpha_i h_i (v_j) = \delta_{ij}$, where
  $\delta_{ij}$ is Kronecker's symbol. In particular, the $v_i$ are
  linearly independent. Hence they form a basis with dual basis
  $\alpha_i h_i$. In a way similar to the computation in
  Equation~\eqref{eq:dual}, we see that
  \[ (\alpha_i h_i \wedge \alpha_j h_j)^\phi = D v_k \wedge v_\ell,
  \]
  where again $D$ is the inverse of $v_0 \wedge v_1 \wedge v_2
  \wedge v_3$ and we choose $k, \ell$ such that $\{i, j, k, \ell\} =
  \{1, 2, 3, 4\}$. Then the voltage of the quadrangle is equal to
  \begin{align*}
    & (v_0 \wedge v_1) (\alpha_0 h_0 \wedge \alpha_1 h_1)
    + (v_1 \wedge v_2) (\alpha_1 h_1 \wedge \alpha_2 h_2) \\
    +\ & (v_2 \wedge v_3) (\alpha_2 h_2 \wedge \alpha_3 h_3)
    + (v_3 \wedge v_0) (\alpha_3 h_3 \wedge \alpha_0 h_0) \\
    =\ & D \bigl((v_0 \wedge v_1) (v_2 \wedge v_3) + (v_1 \wedge
    v_2) (v_0 \wedge v_3) + (v_2 \wedge v_3) (v_0 \wedge v_1) + (v_3
    \wedge v_0) (v_1 \wedge v_2) \bigr) \\
    =\ &0.
  \end{align*}
\end{proof}
\begin{lemma} \label{lm:pentagoninWplusU}
  The voltage of a cycle in $\widetilde H_3 (\mathbb{F})$ of length 5
  is in the $\mathbb{F}_2$-space $W^{(2)} \oplus \langle U \rangle$.
\end{lemma}
\begin{proof}
  Let us take a pentagon $v_0 \tensor h_0, \dotsc, v_4 \tensor h_4$.
  Since for all $i$, we have $v_i \tensor h_i \sim \alpha_i v_i
  \tensor h_i$, we can replace $v_i \tensor h_i$ by $\alpha_i v_i
  \tensor h_i$. This has the effect of subdividing the pentagon into a
  quadrangle (with voltage in $W^{(2)} \oplus \langle U \rangle$) and
  a new pentagon. Now if we prove that the new pentagon has a voltage
  in $W^{(2)} \oplus \langle U \rangle$, then the old one also has
  that property. The process is shown in
  Figure~\ref{fig:pentagon-alpha}. In this way we ensure that $h_i
  (v_i) = 1$ for all $i$.

  \begin{figure}
    \centering
    \myinput{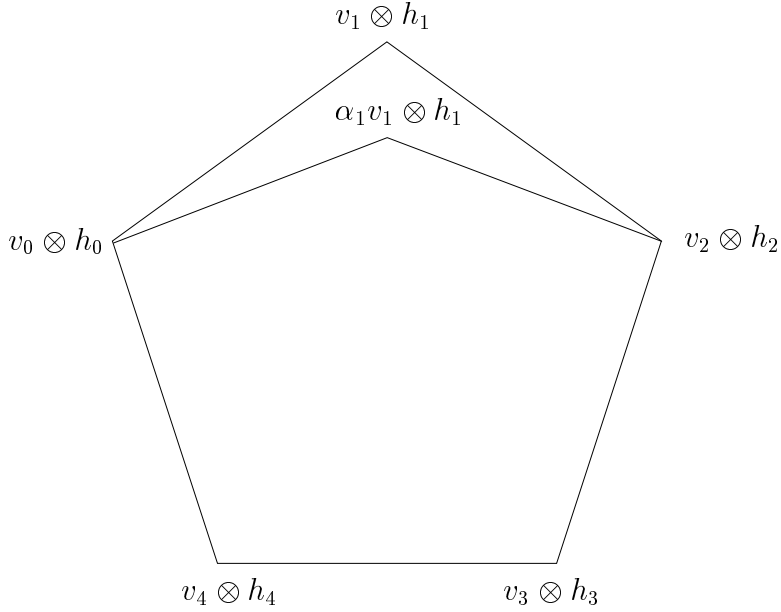}
    \caption{Replacing a vertex in a pentagon by a scalar multiple in
      order to have $h_i (v_i) = 1$.}
    \label{fig:pentagon-alpha}
  \end{figure}
  
  Choose an index $i$ and consider the index set $\{i - 2, i, i + 2\}$
  (modulo $5$). Now suppose that the common null space $\mathcal{N}_i$
  of $h_{i - 2}$, $h_i$ and $h_{i + 2}$ is {\em not} contained in
  $\mathcal{V}_i = \langle v_{i - 2}, v_i, v_{i + 2} \rangle$. Then
  take some $v \in \mathcal{N}_i \setminus \mathcal{V}_i$, and some $h
  \in V^*$ such that $h$ is zero on $\mathcal{V}_i$, but not on $v$.
  Then the vertex $v \tensor h$ is adjacent to $v_{i - 2} \tensor h_{i
    - 2}$, $v_i \tensor h_i$ and $v_{i + 2} \tensor h_{i + 2}$.  Hence
  we can subdivide the pentagon into two quadrangles and a triangle,
  as in Figure~\ref{fig:splitpentagon}.  Therefore its voltage is in
  $W^{(2)} \oplus \langle U \rangle$.

  \begin{figure}
    \centering
    \myinput{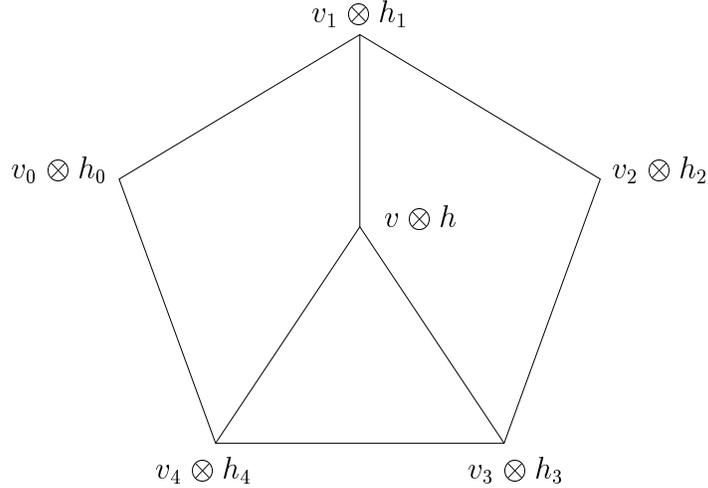}
    \caption{Splitting a pentagon into a triangle and two quadrangles.}
    \label{fig:splitpentagon}
  \end{figure}

  Now suppose that for all indices $i$, we have that $\mathcal{N}_i
  \subseteq \mathcal{V}_i$. Since $\dim \mathcal{N}_i \geq 1$, there
  is a nonzero vector in $\mathcal{V}_i$ on which $h_{i - 2}$, $h_i$
  and $h_{i + 2}$ are all zero. Say this vector is 
  \[ 
  v = \lambda_{i - 2} v_{i - 2} + \lambda_i v_i + \lambda_{i + 2} v_{i + 2}. 
  \]
  If $\lambda_i = 0$, then also 
  \[
  0 = h_{i \pm 2} (\lambda_{i \pm 2} v_{i \pm 2} + \lambda_{i \mp 2}
  v_{i \mp 2}) = \lambda_{i \pm 2},
  \]
  contradicting $v \not= 0$. So we may assume $\lambda_i \not= 0$ and,
  multiplying $v$ by a suitable scalar, $\lambda_i = 1$. Then
  \[
  \begin{array}{r@{\ }c@{\ }l}
    0 = & h_i (v) &= \lambda_{i - 2} h_i (v_{i - 2}) + 1 +
    \lambda_{i + 2} h_i (v_{i + 2}); \\
    0 = & h_{i - 2} (v) & = \lambda_{i - 2} + h_{i - 2} (v_i); \\
    0 = & h_{i + 2} (v) & = \lambda_{i + 2} + h_{i + 2} (v_i).
  \end{array}
  \]
  So we find $\lambda_{i \pm 2} = h_{i \pm 2} (v_i)$ and hence 
  \begin{equation}
    \label{eq:hivis}
    h_{i - 2} (v_i) h_i (v_{i - 2}) + h_{i + 2} (v_i) h_i (v_{i + 2})
    = 1.
  \end{equation}
  If we sum Equation~\eqref{eq:hivis} over all $i$, the right hand
  side is $1$. But every term on the left hand side occurs twice, so
  the left hand side is $0$. Contradiction.
\end{proof}
\begin{corollary} \label{cor:labelsinwplusu}
  The voltage of a cycle in $\widetilde H_3 (\mathbb{F})$ is in the
  $\mathbb{F}_2$-space $W^{(2)} \oplus \langle U \rangle$.
\end{corollary}
\begin{proof}
  By Lemma~$1.3.5$ of Gramlich~\cite{cit:gramlich}, the diameter of
  $H_3 (\mathbb{F})$ is two. It follows that the diameter of
  $\widetilde H_3 (\mathbb{F})$ is also
  two.
  
  Lemmas~\ref{lm:triangleinU},~\ref{lm:quadrangleinWplusU}
  and~\ref{lm:pentagoninWplusU} tell us that the lemma holds for all
  cycles of length at most $5$. Let $c = (v_0 \perp v_1\perp \dotsb
  \perp v_n = v_0)$ be a shortest cycle with a voltage not in $W^{(2)}
  \oplus \langle U \rangle$; so $n \geq 6$.  There is a path of length
  at most $2$ from $v_0$ to $v_3$. Let us call this path $P$. This
  gives us two new cycles: $v_0$, $v_1$, $v_2$, $v_3$, followed by $P$
  in reverse direction, of length at most $5$ and hence with a voltage
  in $W^{(2)} \oplus \langle U \rangle$; and $P$, followed by $v_4,
  \dotsc, v_n$ of length at most $n - 1$ and hence also with a voltage
  in $W^{(2)} \oplus \langle U \rangle$. Therefore the voltage of $c$
  is also in $W^{(2)} \oplus \langle U \rangle$. Contradiction.
\end{proof}
\begin{lemma} \label{lm:allw2occurs}
  For all $w \in W^{(2)}$, there is a cycle in $\widetilde H_3
  (\mathbb{F})$ with voltage $w$.
\end{lemma}
\begin{proof}
  It is sufficient to show that a set of generators of $W^{(2)}$
  occurs as cycles. Note that an $\mathbb{F}$-basis is not necessarily
  sufficient -- we need an $\mathbb{F}_2$-basis.

  We define
  \begin{align*}
    v_1 & = e_1, & v_2 & = e_3, & v_3 & = e_3 + \lambda e_2, \\
    h_1 & = f_3, & h_2 & = f_1, & h_3 & = f_1 + f_4,
  \end{align*}
  and consider the quadrangle $(v_1 \tensor h_2), (v_2 \tensor h_1),
  (v_1 \tensor h_3), (v_3 \tensor h_1)$. Its voltage evaluates to
  $\lambda e_{1,2}^2$. By permuting the basis vectors and by
  choosing different values for $\lambda$, we obtain an
  $\mathbb{F}_2$-basis for $W^{(2)}$.
\end{proof}
\begin{proof}[Proof of Theorem~\ref{thm:mainthm}.]
  We apply Corollary~\ref{cor:subtildelabelingconditions}.
  Lemma~\ref{lm:triangleinU} gives us $T = \langle U
  \rangle_{\mathbb{F}_2}$; then Corollary~\ref{cor:labelsinwplusu} and
  Lemma~\ref{lm:allw2occurs} imply that $M = W^{(2)}$. Hence every
  connected component of $\Gamma$ is a $\card{W^{(2)}}$-fold cover of
  $\Delta$.  Since $W^{(2)} \cong \mathbb{F}^6$, we have finished the
  proof.
\end{proof}

\section{A group of automorphisms} \label{s:autgroup}
Let $\Gamma$ be the graph of Theorem \ref{thm:mainthm} and set $N = S_2 (W) /
\langle U \rangle_{\mathbb{F}_2}$ and $M = W^{(2)}$.  The group $\mathrm{SL}_4 (\mathbb{F})$
acts on $H_3 (\mathbb{F})$, so, by Corollary~\ref{cor:tildelabelledaction},
the group $\mathrm{SL}_4 (\mathbb{F})\ltimes N$ acts on $\Gamma$. By
the results of Section \ref{s:phpgraphs}, an extension of $\mathrm{SL}_4
(\mathbb{F})$ by $M$ acts on a connected component of $\Gamma$.

\begin{theorem}
  Let $\characteristic \mathbb{F} = 2$ and $\card{\mathbb{F}} > 2$.
  Then the group $H$ of
  automorphisms of $\Gamma$ obtained as the stabilizer in
  $\mathrm{SL}_4 (\mathbb{F}) \ltimes N$ of a connected component of $\Gamma$, is a
  non-split extension of $\mathrm{SL}_4 (\mathbb{F})$ by
  $\mathbb{F}^6$.
\end{theorem}
\begin{proof}
  Let $\alpha \in \mathbb{F}$ be an element outside the ground field.
  Let $F$ denote the additive subgroup $\langle 1, \alpha \rangle$ of
  $\mathbb{F}$ of order $4$. For $x \in F$, let $A_x$ be the element
  of $\mathrm{L}_4 (\mathbb{F})$ fixing $e_2$ and $e_4$, and mapping
  $e_1$ and $e_3$ to $e_1 + x e_2$ and $e_3 + x e_4$, respectively.
  We will show that the subgroup $A = \{ A_x
  \mid x \in F\}$ does not lift to a subgroup of $H$.
  
  Let $\widetilde v_x = (e_1 + x e_2, f_1) \in V(\widetilde \Delta)$
  and let $v_x$ be the corresponding vertex $(\langle e_1 + x
  e_2\rangle, \langle f_1 \rangle) \in V (\Delta)$; similarly, let
  $\widetilde u = (e_3, f_3) \in V(\widetilde \Delta)$ and let $u$ be
  the corresponding vertex $(\langle e_3 \rangle, \langle f_3 \rangle)
  \in V(\Delta)$. Then $v_x \perp u$ for all $x \in F$.

  We define a basis for $W$, and write down the images under $A_x$.
  \[
  \begin{array}{*{6}{c}r}
    w_1 & = & e_1 \wedge e_2 & = & (f_3 \wedge f_4)^\phi & \mapsto & w_1; \\
    w_2 & = & e_1 \wedge e_3 & = & (f_2 \wedge f_4)^\phi & \mapsto &
    w_2 + x w_3 + x w_4 + x^2 w_5; \\
    w_3 & = & e_1 \wedge e_4 & = & (f_2 \wedge f_3)^\phi & \mapsto &
    w_3 + x w_5; \\
    w_4 & = & e_2 \wedge e_3 & = & (f_1 \wedge f_4)^\phi & \mapsto &
    w_4 + x w_5; \\
    w_5 & = & e_2 \wedge e_4 & = & (f_1 \wedge f_3)^\phi & \mapsto & w_5; \\
    w_6 & = & e_3 \wedge e_4 & = & (f_1 \wedge f_2)^\phi & \mapsto & w_6.
  \end{array}
  \]
  The voltage along the dart between $u$ and $v_x$ is
  \[ 
  \bigl((e_1 + x e_2) \wedge e_3 \bigr)(f_1 \wedge f_3) = w_2 w_5 + x
  w_4 w_5.
  \]
  In the terminology of Lemma~\ref{lm:GM} we take $v = v_0$, so $v^{A_x} = v_x$. For all $x \in F$, we choose
  the path from $v_x$ to $v$ to run through $u$. Then $\lambda (A_x) =
  x w_4 w_5$. Hence 
  \[
  f (A_x, A_y) = \lambda (A_{x + y}) + \lambda (A_x)^{A_y} + \lambda
  (A_y) = x w_4 w_5 + x w_4^{A_y} w_5^{A_y} = x w_4 w_5 + x (w_4 + y
  w_5) w_5 = x y w_5^2.
  \]
  So we can construct the extension of $A$ as a group defined on the
  set $F \times M$, where $[x, m]$ stands for the element $(A_x,
  \lambda (A_x) + m)$ in $\mathrm{SL}_4 (\mathbb{F}) \ltimes N$; the
  multiplication rule is then $[x, m] [y, n] = [x + y, x y w_5^2 +
  m^{A_y} + n]$. Now suppose that $A$ lifts to a subgroup of $H$, that is, there is a 
  function $c\colon F \to M$ such that $\{[x, c (x)] \mid x \in F\}$
  is a group isomorphic to $A$. Then $[1, c (1)]$, $[\alpha, c
  (\alpha)]$, and $[\alpha + 1, c (\alpha + 1)]$ need to have order $2$.
  
  Now $[x, m]^2 = [0, x^2 w_5^2 + m^{A_x} + m]$, so $[x, m]$ has order
  two if and only if $m^{A_x} + m = x^2 w_5^2$. By elementary linear
  algebra we find that this is true for $x \not= 0$ if and only if $m
  \in w_3^2 + S$, where $S = \langle w_1^2, w_3^2 + w_4^2, w_5^2,
  w_6^2 \rangle_{\mathbb{F}}$. Note that $S$ is $A$-invariant.
  
  Let $s (x) = w_3^2 + c (x)$. Then $s (x) \in S$. Since $[1, c (1)]
  [\alpha, c (\alpha)] = [\alpha + 1, c (\alpha + 1)]$, we have $s
  (\alpha + 1) = w_3^2 + \alpha w_5^2 + c (1)^{A_\alpha} + c (\alpha)
  = w_3^2 + (\alpha + \alpha^2) w_5^2 + s (1) + s (\alpha) \not\in
  S$, a contradiction.
  
  Since $A$ does not lift to a subgroup of $H$, neither does $\mathrm{SL}_4
  (\mathbb{F})$. In other words, the extension of $\mathrm{SL}_4 (\mathbb{F})$
  by $M$ is non-split.

\end{proof}

\end{document}